\documentclass[12pt]{article}
\usepackage{fullpage}
\usepackage{times}
\usepackage{amsmath}
\usepackage{latexsym}
\usepackage{amssymb,amsfonts,amsthm}
\usepackage{graphicx}
\usepackage{graphics}
\usepackage{float}
\usepackage{multirow}
\usepackage{multicol}
\usepackage{rotating}
\usepackage{setspace}
\usepackage{mathrsfs}
\usepackage{tikz}
\usepackage{tabularx}
\usepackage{listings}
\usepackage{hyperref}
\hypersetup{
    colorlinks=true,
    linkcolor=blue,
    filecolor=magenta,      
    urlcolor=red,
}
\usepackage[utf8]{inputenc}
\usepackage{wrapfig}
\usepackage{parskip}
\usetikzlibrary{positioning,chains,fit,shapes,calc}

 \renewcommand{\-}{\ensuremath{\text{--}}}

\let\oldmarginpar\marginpar
\renewcommand\marginpar[1]{\-\oldmarginpar[\raggedleft\footnotesize #1]%
{\raggedright\footnotesize #1}}

\newtheorem{theorem}{Theorem}

\begin{document}

\title{On Projective Planes of Order 16 Associated with 1-rotational 2-$(52,4,1)$ Designs}

\author{Zazil Santizo Huerta, Vladimir Tonchev, and Melissa Keranen\\
Michigan Technological University\\
Department of Mathematical Sciences\\
Houghton, MI 49931, USA\\
}

\maketitle

 \begin{abstract}
A maximal arc of degree $k$ in a finite projective plane $\cal{P}$ of order $q=ks$
is a set of $(q - s +1)k$ points that meets every line of $\cal{P}$ in either $k$ or $0$ points. The collection of the nonempty intersections of a maximal arc 
with the lines of $\cal{P}$  is a resolvable Steiner  2-$((q - s +1)k, k, 1)$ design. 
Necessary and sufficient conditions for a resolvable Steiner 2-design to be embeddable as a maximal arc in a projective plane were proved recently in \cite{T}. 
Steiner designs associated with  maximal arcs in the known projective planes of order 16 were analyzed in \cite{GWT}, where it was shown that some of the associated designs are embeddable in two nonisomorphic planes. 

Using MAGMA, we conducted an analysis to ascertain whether any of the 22 non-isomorphic 1-rotational 2-$(52,4,1)$ designs, previously classified in \cite{B}, could be embedded in maximal arcs of degree 4 within projective planes of order 16. This paper presents a summary of our findings, revealing that precisely only one out of the the twenty-two 1-rotational designs from \cite{B} is embeddable in a plane of order 16, being the Desarguesian plane $PG(2,16)$. 

 \end{abstract}

\section{Introduction}

We assume familiarity with basic facts and notions from combinatorial design
theory, and emphasize the following definitions needed to state and prove our main theorem. 

A 2-$(v,k,\lambda)$ design $D=\{ X, \cal{B} \}$
is a pair of a {\bf point set} $X$ of size $v$ and a collection of $k$-subsets $\cal{B}$=$\{ B_j \}_{j=1}^b$ called {\bf blocks}, such that every pair of points is contained in exactly $\lambda$ blocks. Every point of $X$ is contained in $r=\lambda(v-1)/(k-1)$ blocks, and the number of blocks is $b={\frac{v(v-1)}{k(k-1)}}\lambda$. A design with $v=b$ is called symmetric, and a Steiner design is a design with $\lambda=1$.

Let $D=\{X,B\}$ be a design where $X = \{x_1, \ldots, x_v \}$ and $B = \{B_1, \ldots, B_b \}$. The incidence matrix of $D$ is the $v \times b$ 0-1 matrix $M =(m_{i,j})$ defined by the rule
\[
m_{i,j}= \begin{cases} 
      1 & \hbox{if} \ x_i \in A_j\\
      0 & \hbox{if} \ x_i \notin A_j 
      \end{cases}
\]

The incidence matrix, $M$, of a $(v, b, r, k, \lambda)$-BIBD satisfies the following
properties:

\begin{enumerate}
\item every column of $M$ contains exactly $k$ “1”s;
\item every row of $M$ contains exactly $r$ “1”s;
\item two distinct rows of $M$ both contain “1”s in exactly $\lambda$ columns.
\end{enumerate}

Two designs are {\it isomorphic} if there is a bijection between their point sets that maps every block of the first design to a block of the second design. An {\it automorphism} of a design is any isomorphism of the design to itself. The set of all automorphisms of $D$ form the automorphism group $Aut(D)$ of $D$.

A parallel class of a 2-$(v,k,\lambda)$ design $D$ with $v=nk$ is a set of $n=v/k$ pairwise disjoint blocks. A resolution of $D$ is a partition of the collection of blocks
$\cal{B}$ into $r=\lambda(v-1)/(k-1)$ disjoint  parallel classes. A design is  resolvable if it admits at least one resolution.

\

A finite projective plane of order $q$ is defined as a set of $q^2+q+1$ points, a set of lines, and a relation between points and lines, having the following properties:

\begin{itemize}

  \item Any two points determine a line,
  \item Any two lines determine a point,
  \item Every point has $q+1$ lines through it, and
  \item Every line contains $q+1$ points.

\end{itemize}

 A projective plane of order $q\ge 2$ is a symmetric Steiner 2-$(q^2 + q + 1,q+1,1)$ design. An affine plane of order $q\ge 2$ is a Steiner 2-$(q^2,q,1)$ design.

A maximal $(m,k)$-arc in a projective plane  $\cal{P}$ of order $q=ks$ is a set of $m=(q - s +1)k$ points that meets every line of $\cal{P}$ in either $k$ or $0$ points. The set of lines of $\cal{P}$ which have no points in a given maximal $((sk-s+1)k, k)$-arc   $\cal{A}$ determines a {\bf dual} maximal  $((sk-k+1)s, s)$-arc ${\cal{A}}^{\perp}$ in the dual plane ${\cal{P}}^{\perp}$.

The classical (or Desarguesian) plane $PG(2, p^t)$ of order $n=p^t$, where $p$ is prime and $t \geq 1$,
has as points the 1-dimensional subspaces of the 3-dimensional vector space $V^3$ over the finite field of order $p^t$, and as blocks (or lines), the 2-dimensional subspaces of $V^3$.

\

In Section 2, we provide a framework for understanding the relationships between maximal arcs and Steiner 2-designs, resolutions, and their properties within projective planes. It establishes theoretical groundwork and presents specific findings within the context of planes of order 16.

In Section 3, we give the construction and analyze the 22 non isomorphic resolvable 1-rotational 2-$(52, 4, 1)$ designs to determine whether it is possible to embed them as maximal arcs of degree 4 within projective planes of order 16.

In section 4, we provide a clear outline of the computational approach employed to determine the embeddability of the given 1-rotational 2-$(52,4,1)$ designs into projective planes of order 16. The analysis of our computations culminates in the presentation of the main theorem, which identifies exactly one design is embeddable. 

The following theorems concern the existence of maximal arcs within Desarguesian planes.

\begin{theorem} \cite{D} The Desarguesian plane $PG(2,2^m)$ contains  maximal arcs for every
$k=2^i$, $1\le i < m$.
\end{theorem}

\begin{theorem}\cite{BBM} Maximal $k$-arcs with $1< k < q$ do not exist in any Desarguesian plane, $PG(2,q)$, if $q$ is odd.
\end{theorem}

\section{Previous Results}

Let $\cal{P}$ be a projective plane of order $q=sk$, $1< k < q$.

\begin{itemize}
\item The nonempty intersections of a maximal $((sk-s+1)k, k)$-arc
$\cal{A}$
with the lines of $\cal{P}$ form a Steiner 2-$((sk-s+1)k, k, 1)$ design $D$.
\item Respectively, the dual  $((sk-s+1)k, k)$-arc   $\cal{A}^\perp$
is the point set of a Steiner 2-$((sk-k+1)s, s, 1)$ design $D^{\perp}$.
\end{itemize}

\begin{theorem}
\label{3}
\

\begin{itemize}
\item The points of  ${\cal{A}}^{\perp}$ determine a set of $(sk-k+1)s$
resolutions of $D$.
\item The points of ${\cal{A}}$ determine a set of $(sk-s+1)k$
resolutions of the design $D^{\perp}$ associated
with the dual arc ${\cal{A}}^{\perp}$.

\end{itemize}
\end{theorem}

Every two  resolutions of $D$ (resp.  $D^{\perp}$)
share one parallel class.

\

Tonchev gave the following definition in \cite{GWT}, which is motivated by properties of  Steiner designs associated with maximal arcs.

Two resolutions, $R_1$, $R_2$ of a  2-$((sk -s +1)k, k, 1)$ design,
\[ R_1 = P^{(1)}_1 \cup P^{(1)}_2 \cup \cdots P^{(1)}_r, \  R_2 =  
 P^{(2)}_1 \cup P^{(2)}_2 \cup \cdots P^{(2)}_r 
\]
are {\bf compatible} if they share one parallel class, $P^{(1)}_i = P^{(2)}_j$, and  $|P^{(1)}_{i'} \cap P^{(2)}_{j'}|\le 1$ 
for  $(i',j') \neq (i,j)$. More generally, a {\bf set of $m$ resolutions
$R_{1}, \ldots, R_{m}$ is compatible} if every two of these resolutions are  compatible.

A Steiner 2-$((sk-s+1)k, k, 1)$ design associated with a maximal $((sk-s+1)k, k)$-arc  $\cal{A}$ in a plane $\cal{P}$ of order $q=ks$ admits a set of  $(sk-k+1)s$ pairwise compatible resolutions.

\

The following result tells us  when a design is embeddable in a projective plane as a maximal arc.

\begin{theorem} \cite{T} 
\label{CR} Let $S=\{ R_1,\ldots, R_m \}$ be a set of $m$ pairwise  compatible resolutions of a  2-$((sk -s +1)k, k, 1)$ design $D=\{ X, \cal{B} \}$. Then
\[ m\le (sk-k+1)s. \]
The equality $m=(sk-k+1)s$ holds if and only if there exists a projective plane $\cal{P}$ of order $q=sk$ such that $D$ is embeddable in $\cal{P}$ as a maximal $((sk-s+1)k, k)$-arc.
\end{theorem}

Hamilton, et. al. have determined which of the known projective planes of order 16 contain maximal $(52,4)$-arcs.

\begin{theorem} \cite{G} All but possibly four of the 22 known projective planes of order 16, (BBH2, BBS4, and their duals), contain maximal $(52,4)$-arcs.
\end{theorem}

\

The following theorem states that in the projective geometry over the Galois field $GF(2)$ with 16 elements, there are exactly two distinct sets of $(52, 4)$-arcs. These arcs are {\it projectively inequivalent}, meaning they cannot be transformed into each other through projective transformations while preserving the structure. Furthermore, both of these arcs belong to a specific type known as {\it Denniston} type.

\begin{theorem} \cite{Ball} $PG(2,16)$ contains exactly two projectively inequivalent maximal $(52,4)$-arcs, both of Denniston type.
\end{theorem}

In the following theorem, the listed planes contain $(52,4)$ maximal arcs, which have associated 2-$(52, 4, 1)$ designs. Importantly, each design admits two different sets of 52 compatible resolutions, leading to two non-isomorphic planes of order 16.

\begin{theorem} \cite{GWT} The planes $LMRH$, JOHN, BBH1, and JOWK, contain maximal $(52,4)$ arcs whose associated 2-$(52,4,1)$ designs
admit {\bf two different } sets of 52 compatible resolutions
that lead to {\bf two nonisomorphic} planes of order 16.

\end{theorem}

\

We now study a different question. Instead of asking which projective planes contain a maximal arc, we would like to answer whether a given design can be embedded as a maximal arc in a projective plane.

\section{1-Rotational 2-$(52,4,1)$ Designs}

A 2-$(v,k,\lambda)$ design is 1-rotational if it admits an automorphism of order $v-1$.

Buratti and Zuanni \cite{BZ} have shown that a 1-rotational $(52,4,1)$-RBIBD is equivalent to a 4-set $F=\{A_1, A_2, A_3, A_4\}$ of 4-subsets of $\mathbb{Z}_{51}$ satisfying that any element of $\mathbb{Z}_{51}$ - $\{0,17,34\}$ is representable as a difference of two elements of a same $A_i$, and $A_1 \cup A_2 \cup A_3 \cup A_4 = \mathbb{Z}_{17} - \{0\} \pmod{17}$. They also provided a construction of the 22 non-isomorphic 1-rotational 2-$(52,4,1)$ designs in \cite{B}.

\begin{theorem}\cite{B}
Up to isomorphism, there exist exactly 22 nonisomorphic
resolvable 1-rotational 2-$(52,4,1)$ designs.
\end{theorem}

Given a $(51,3,4,1)$ resolvable difference family $F = \{A_1, A_2, A_3, A_4\}$. We give the associated 2-$(52,4,1)$-RBIBD. This design has the point set $\mathbb{Z}_{51} \cup \{\infty \}$ and parallel classes obtained by developing $\pmod{51}$ the starter parallel class given by 

$$\{0,17,34,\infty\} \cup \{A_i+17j | 1 \leq i \leq 4; 0 \leq j<3 \}$$

A complete set of representatives $F_i$, $i=1, \ldots, 22$, for each of the 22 1-rotational designs are listed in \cite{B}.
\

We would like to determine which of the Buratti-Zuanni 2-$(52,4,1)$ designs are embeddable as maximal $(52,4)$-arcs in a projective plane of order 16.

\section{Main Theorem and Results}

We wish to determine which of these designs can be embedded into a projective plane of order 16 as a maximal $(54,4)$- arc. Thus, to answer our question, we need that the number of compatible resolutions in the design will be 52, by Theorem \ref{CR}.

Finding the compatible sets of maximum size, each containing 52 resolutions, is essentially equivalent to identifying subsets of resolutions where every pair is compatible. We now present some fundamental definitions of graph theory, which are essential for establishing the graphical theoretical equivalent of this problem. 

A {\bf graph} is a pair $G = (V, E)$ (sometimes called an undirected graph) consists of a set $V$ of vertices (also called nodes) and a set $E$ of edges. If two vertices in a graph are connected by an edge, we say the vertices are {\bf adjacent}. A {\bf clique} is a subset of vertices of an undirected graph such that every two distinct vertices in the clique are adjacent. By constructing a graph $P$ where each vertex represents a resolution of $F_i$, and two vertices are adjacent if the corresponding resolutions are compatible, we effectively represent the compatibility relationships between resolutions as edges in the graph. Since a clique in a graph is a subset of vertices where each pair is adjacent, finding maximum-sized cliques in graph $P$ corresponds to identifying sets of maximum size among resolutions that are compatible. This is because each clique represents a compatible set of resolutions, and its size determines the number of resolutions in the set. Therefore, identifying the size of the largest cliques in $P$ concurrently reveals the size of the largest set of compatible resolutions.

\

We now proceed to describe the outline of the algorithm created to compute this number.

\

{\bf Algorithm Outline:}
\begin{itemize}
\item Construct the incidence matrices  of the twenty-two  1-rotational designs from \cite{B}.  
 \item Using Magma, compute all parallel classes and resolutions of each design $F_i$. 
\item For each $F_i$, define a graph $P$ whose vertices are the resolutions of  $F_i$, where two resolutions are adjacent if and only if they are compatible. 
\item Determine the largest size, $m$, of a clique in $P$.
\item If $m=52$, then by Theorem \ref{CR}, $F_i$ is embbedable as a maximal $(52,4)$-arc in the projective plane.
\end{itemize}

\

Table \ref{T1} summarizes the results of running the magma code provided in section 5. We see that $F_{17}$ is the only design that achieves the maximum clique size of 52. This leads us to our main theorem.

\begin{theorem}

Only one of the twenty-two  1-rotational 2-(52,4,1) designs,
namely $F_{17}$, is embeddable as a maximal (52,4)-arc in a plane of order 16, being the Desarguesian plane $PG(2,16$).
\end{theorem}

\begin{table}
\centering
\begin{tabularx}{0.8\textwidth} { 
  | >{\centering\arraybackslash}X 
  | >{\centering\arraybackslash}X 
  | >{\centering\arraybackslash}X
  | >{\centering\arraybackslash}X
  | >{\centering\arraybackslash}X | }
  \hline
{\bf  Design\cite{B}} & ${\bf |Aut(F_i)|}$ & {\bf Par. Cl.} & {\bf Res.} & {\bf Comp. Res.} \\
  \hline
 $F_1$ & 51 & 34 & 1 & 1 \\
  \hline
 $F_2$ & 51 & 17 & 1 & 1 \\
  \hline
 $F_3$ & 51 & 17 & 1 & 1 \\
  \hline
 $F_4$ & 51 & 17 & 1 & 1 \\
  \hline
 $F_5$ & 51 & 34 & 1 & 1 \\
  \hline
 $F_6$ & 51 & 17 & 1 & 1 \\
  \hline
 $F_7$ & 51 & 34 & 1 & 1 \\
  \hline
 $F_8$ & 51 & 17 & 1 & 1 \\
  \hline
 $F_9$ & 51 & 34 & 1 & 1 \\
  \hline
 $F_{10}$ & 51 & 34 & 1 & 1 \\
  \hline
 $F_{11}$ & 51 & 17 & 1 & 1 \\ 
  \hline
$F_{12}$ & 204 & 442 & 18 & 1 \\
  \hline
 $F_{13}$ & 204 & 289 & 18 & 1 \\
  \hline
 $F_{14}$ & 204 & 153 & 1 & 1 \\
  \hline
 $F_{15}$ & 102 & 306 & 1 & 1 \\
  \hline
 $F_{16}$ & 408 & 51  & 1 & 1 \\
  \hline
 $F_{17}$ & 408 & 2550 & 460 & 52 \\
  \hline
 $F_{18}$ & 51 & 34 & 1 & 1 \\
  \hline
 $F_{19}$ & 51 & 34 & 1 & 1 \\
  \hline
 $F_{20}$ & 51 & 17 & 1 & 1 \\
  \hline
 $F_{21}$ & 51 & 17 & 1 & 1 \\
  \hline
 $F_{22}$ & 51 & 17 & 1 & 1 \\
  \hline
\end{tabularx}
\caption{The 22 Resolvable 1-rotational 2-$(52,4,1)$ designs}
\label{T1}
\end{table}

\newpage
\section{Magma Code}

\begin{lstlisting}
> DifferenceSet := func< S, n| [{(x + i) mod n + 1 : x in S } : i in [1 .. n] ] >;      
> A := {18,33,22,46};
> B := {21,30,37,31};
> C := {6,45,25,43};
> D := {24,27,19,49};
> E:= { {1,18,35,52}, {2,19,36,52}, {3,20,37,52}, {4,21,38,52}, {5,22,39,52}, {6,23,40,52}, {7,24,41,52}, {8,25,42,52}, {9,26,43,52}, {10,27,44,52}, {11,28,45,52}, {12,29,46,\
52}, {13,30,47,52}, {14,31,48,52}, {15,32,49,52}, {16,33,50,52}, {17,34,51,52} };
> F := Design< 2, 52| DifferenceSet(A, 51), DifferenceSet(B, 51), DifferenceSet(C, 51), DifferenceSet(D, 51), E >;
> pRank(F, 2);
41
> AutomorphismGroup(F);
Permutation group acting on a set of cardinality 52
Order = 408 = 2^3 * 3 * 17
    (1, 50, 49, 23, 10, 29, 13, 5)(2, 25, 11, 4, 26, 37, 17, 7)(3, 51, 24, 36, 42, 45, 21, 9)(6, 27, 12, 30, 39, 18, 33, 15)(8, 28, 38, 43, 20, 34, 41, 19)(14, 31)(16, 32, 40, 
        44, 46, 47, 22, 35)
    (2, 27, 14, 33, 17, 9, 5, 3)(4, 28, 40, 46, 49, 25, 13, 7)(6, 29, 15, 8, 30, 41, 21, 11)(10, 31, 16, 34, 43, 22, 37, 19)(12, 32, 42, 47, 24, 38, 45, 23)(18, 35)(20, 36, 44,
        48, 50, 51, 26, 39)
> G :=AllParallelClasses(F);
> #G;
2550
> H := AllParallelisms(F);
> #H;
460
> V := { H[i]: i in [1 .. #H] };
> edgesSet := {};
> for i in [1 .. #H] do
for> for j in [1 .. #H] do
for|for> T := {};
for|for> for m in [1 .. 17] do
for|for|for> for n in [1 .. 17] do
for|for|for|for> if (i lt j) and (#(H[i] meet H[j]) eq 1) and ((#((SetToIndexedSet(H[i])[m]) meet (SetToIndexedSet(H[j])[n])) le 1)) then
for|for|for|for|if> T join:= {<m,n>};
for|for|for|for|if> end if;
for|for|for|for> end for;
for|for|for> end for;
for|for> if #T eq 288 then
for|for|if> edgesSet join:= {{i, j}};
for|for|if> end if;
for|for> end for;
for> end for;
> #edgesSet;
1326
> P := Graph< V| edgesSet >;
> M := MaximumClique(P);
> W := AllCliques(P, #M);
> #M;
52
> #W;
1

\end{lstlisting}

\section{Conclusions and Future Work}

Our goal is to validate our magma code by applying it to other 2-$(v,k,\lambda)$ 1-rotational designs that meet the criteria for embedding as a maximal arc in a projective plane. Additionally, we aim to extend the applicability of this code to test it on other types of Steiner 2-designs.


\begin{thebibliography}{99}

\bibitem{Ball} Ball, S. M., and Blokhuis, A. (2002). The classification of maximal arcs in small Desarguesian planes. Bulletin of the Belgian Mathematical Society : Simon Stevin, 9(3), 433-445.

\bibitem{BBM} S. Ball, A. Blokhuis, F. Mazzocca, Maximal arcs in Desarguesian planes of odd order to not exist, Combinatorica 17 (1997) 31-41.

\bibitem{B} M. Buratti and F. Zuanni, The 1-rotational $(52,4,1)$-RBIBD's,
{\it JCMCC} {\bf 30} (1999), 99 - 102. 

\bibitem{BZ} Buratti, Marco and Zuanni, Fulvio. (2000). Addendum to: "{$G$}-invariantly resolvable Steiner 2-designs which are 1-rotational over {$G$}". Bulletin of the Belgian Mathematical Society, Simon Stevin. 5. 10.36045/bbms/1103055695. 

\bibitem{D} R.H.F. Denniston, Some maximal arcs in finite projective planes, J. Combin. Theory 6 (1969) 317-319.

\bibitem{GWT} M. Gezek, T. Wagner and V. D. Tonchev,
Maximal arcs in projective planes of order 16 and related designs,
{\it Advances in Geometry}, to appear.

\bibitem{G} Gezek, Mustafa, "COMBINATORIAL PROBLEMS RELATED TO CODES, DESIGNS AND FINITE GEOMETRIES", Campus Access Dissertation, Michigan Technological University, 2017.

\bibitem{T} V. D. Tonchev, On resolvable Steiner 2-designs and maximal arcs in projective planes, {\it Designs, Codes and Cryptography}, {\bf 84},
July (2017), 165-172.




\end{thebibliography}
\end{document}